\documentclass[12pt,twoside]{myarticle}                                   

\usepackage{theorem}
\usepackage{epsfig}
\usepackage{amsfonts}
\usepackage{amssymb}

\theoremstyle{change}{\theorembodyfont{\slshape}
\newtheorem{theorem}{Theorem.}[section]
\newtheorem{proposition}[theorem]{Proposition.}
\newtheorem{lemma}[theorem]{Lemma.}

\newtheorem{corollary}[theorem]{Corollary.}

}

\theoremstyle{change}{
\theorembodyfont{\rmfamily}

}

\def\proof{\noindent{\bf Proof.}\enspace}
\def\endproof{ \quad $\kasten$\bigskip}

\def\F{\mathop{\cal{F}{}}}
\def\G{\mathop{\cal{G}{}}}
\def\I{\mathop{\cal{J}}\nolimits}
\def\H{\mathop{\cal{H}{}}}

\def\A{\mathop{\cal{A}^{\bullet}}}
\def\E{\mathop{\cal{E}^{\bullet}}}
\def\randsigma{\partial\sigma}
\def\randtau{\partial\tau}
\def\m{\mathfrak{m}}
\def\IH{\mathop{\hbox{$I$}\!\hbox{$H$}}\nolimits}
\def\Pol{P}
\def\IP{u}
\def\IQ{v}

\def\epsilon{\varepsilon}

\def\b#1{\overline{#1}}

\def\ZZ{{\mathbb Z}}
\def\RR{{\mathbb R}}

\def\mal{\mathbin{\! \cdot \!}}

\def\cone{\mathop{\hbox{\rm cone}}}

\def\id{\mathop{\rm id}\nolimits}

\def\GL{\mathop{\rm GL}\nolimits}
\def\osubset{\subset \kern-10pt {\rm o}\ } 

\def\kasten{\mathord{\vbox{\hrule
                     \hbox{\vrule
                     \hskip5pt
                     \vrule height5pt
                     \vrule}
                     \hrule}}}
\def\text#1{\hbox{\rm #1}}

\def\bigtopmapright#1{\smash{\mathop{\hbox to 35pt{\rightarrowfill}}\limits^{#1}}}
\def\rmapdown#1{\Big\downarrow\rlap{$\vcenter {\hbox{$\scriptstyle#1$}}$}}
\def\lmapdown#1{\llap{$\vcenter{\hbox {$\scriptstyle#1$}}$}\Big\downarrow}

\begin{document}

\title{Lower Bounds for the Generalized $h$-Vectors of Centrally Symmetric Polytopes}
\author{Annette A'Campo--Neuen%
\footnote{\small e-mail: acampo@mathematik.uni-mainz.de}
}

\maketitle

\pagestyle{myheadings}
\markboth{\hfill A.~A'Campo--Neuen\hfill}{\hfill Lower Bounds for the Generalized $h$-Vectors of  Centrally Symmetric Polytopes \hfill}

\begin{abstract} 
\noindent In a previous article, we proved
 tight lower bounds for the coefficients of
the generalized $h$-vector of a centrally symmetric rational polytope using intersection cohomology of the associated projective toric variety.
Here we present a new proof based on the theory of combinatorial intersection cohomology developed by Barthel, Brasselet, Fieseler and Kaup. This theory is also valid for nonrational polytopes when there are no corresponding toric varieties. So we can establish  our bounds for centrally symmetric polytopes even without requiring them to be rational. 
\end{abstract}

\section*{Introduction}
In \cite{St1}, R.~Stanley proved tight lower bounds for the $h$-vector of 
a simplicial centrally symmetric polytope. The entries of the $h$-vector
are linear combinations of face numbers of the polytope, and they determine
the face numbers completely.

Stanley also introduced the {\it generalized $h$-vector\/} of an arbitrary
convex polytope (see \cite{St2}). It is  a combinatorial
invariant of the polytope defined by recursion over its faces, 
and in the simplicial case it coincides with the usual
 $h$-vector.

If the polytope $\Pol$ is
rational then there is an associated projective toric variety $X_{\Pol}$,
and the coefficients of the generalized $h$-vector of ${\Pol}$ have a
topological interpretation as
Betti numbers of the intersection cohomology of $X_{\Pol}$. Let $n$ denote
the dimension of the polytope. It follows from   Poincar\'e-duality that the generalized $h$-vector $(h_0,\dots,h_n)$ is palindromic (i.e. $h_j=h_{n-j}$ for all $j$). The Hard Lefschetz Theorem for
intersection cohomology implies that the $h$-vector  is unimodal (i.e.
the coefficients  increase  up to the middle coefficient(s) and then decrease again).

In the article \cite{hVektoren}, we considered among others centrally symmetric rational polytopes and asked for combinatorial conditions imposed on the generalized $h$-vector by the existence of the central symmetry. Using the topological
interpretation  via intersection cohomology,
we proved tight lower bounds for the coefficients of the generalized $h$-vector of a centrally symmetric rational polytope. In the simplicial case these are precisely the bounds of Stanley (see \cite{St1}). 

The aim of this article is to show that the same bounds remain valid even
if we do not assume the centrally symmetric polytope to be rational. 
Our proof is based on the theory of combinatorial intersection cohomology
of fans developed by  G.~Barthel, J.-P.~Brasselet, K.~Fieseler and
L.~Kaup (see \cite{virtual}) and independently by P.~Bressler and Lunts (see
\cite{Bressler}).


They discovered  that one can completely characterize the 
 intersection cohomology
of a toric variety by combinatorial and algebraic data associated
to the corresponding fan, namely in terms of a {\it minimal extension
sheaf\/} on the fan considered as a topological space where the subfans
are the open subsets (see \cite{Festschrift}). Associating
 an analogous object to a non--rational fan, they define
a combinatorial intersection cohomology satisfying
similar formal properties as the usual intersection cohomology.

Barthel, Brasselet, Fieseler and Kaup conjectured that for the
combinatorial intersection cohomology of a fan arising from a polytope,  
 a combinatorial version of the Hard Lefschetz theorem holds. Moreover they proved that if such a Hard Lefschetz theorem is  true  then the odd Betti numbers of the combinatorial intersection cohomology of a polytopal fan vanish and the even Betti-numbers are precisely 
 the coefficients of the generalized $h$-vector of the corresponding polytope (see \cite{virtual}).
 
The Hard Lefschetz theorem in this context was recently proved by Kalle Karu (see \cite{Karu}). The fact that a combinatorial Hard Lefschetz theorem holds has striking consequences. For example, the coefficients of generalized $h$-vector of an arbitrary polytope are non-negative which is not at all clear from the definition. Moreover, the generalized $h$-vector of an arbitrary polytope is  unimodal.

We apply these results to a centrally symmetric polytope $\Pol$ of dimension $n$. Denoting its generalized $h$-vector by $(h_0,\dots,h_n)$, we prove the following for the polynomial $h_P:=\sum_{j=0}^n h_jx^j$ (see Theorem \ref{Hauptsatz}):

\noindent{\sc Theorem.\enspace} {\sl 
If a polytope $\Pol$ of dimension $n$
admits a central
symmetry then the polynomial
$$h_{\Pol}(x)-(1+x)^n$$
has nonnegative, even coefficients, it is palindromic and unimodal.
That means that we have
the following bounds for the coefficients
$h_j$ of $h_{\Pol}$:
$$h_j-h_{j-1}\geq {n\choose j}- {n\choose j-1}\quad\text{for $j=1,\dots,n/2$.}
$$}

\bigskip

Note that $(1+x)^n$ occurs as the $h$-polynomial of the $n$-dimensional cross-polytope. We can reformulate the lower bounds given by the theorem in terms of the partial ordering on real polynomials of degree $n$ defined by  coefficientwise comparison, i.e. $a=\sum_{j=0}^n a_j x^j \leq b=\sum_{j=0}^n b_j x^j$ if and only if $a_j\leq b_j$ for all $j$.  
The $h$-polynomial of the $n$-dimensional cross-polytope is minimal in this sense and in fact this is the only polytope realizing the minimum (see Corollary \ref{Minimum}).

\noindent{\sc Corollary.\enspace} {\sl
Let $\Pol$ be an $n$-dimensional centrally symmetric polytope. Then
$$h_{\Pol}\geq (1+x)^n\,.$$ 
Moreover, equality holds if and only if the polytope $\Pol$ is affinely equivalent to the $n$-dimensional cross-polytope.}

\bigskip

\section{Preliminaries}

Let $\Pol$ denote a convex polytope of dimension $n$ in  $V:=\RR^n$.
Assume that zero lies in the interior of $\Pol$. Then the polytope $\Pol$
defines a complete fan in $V$ consisting of the
cones through its proper faces:
$$\Delta_{\Pol}:=\{\cone(F); F\text{ proper face of } \Pol\} \cup \{0\}\,,$$
Moreover, this fan is equipped with a {\it strictly concave support
function\/}, that means a concave function whose restriction to any
cone in $\Delta_{\Pol}$ is linear and such that for any two different maximal
cones the linear functions obtained by restriction are different.
To define this function, consider
the dual polytope 
${\Pol}^*:=\{u\in V^* ; \langle u,v\rangle \leq -1 \text{ for all } 
v\in {\Pol}\}$
of $\Pol$.
There is an order--reversing one--to--one--correspondence
between the proper faces of $\Pol$ and the
proper faces of $\Pol^*$ defined by
$$F\mapsto s_F:=\{u\in \Pol^* ; \langle u,v\rangle=-1 \text{ for all } v\in F\}
\,.$$
So the vertices of $\Pol^*$ are of the form $s_F$, where $F$ is
a one--codimensional face of $\Pol$, and $s_F$ defines a linear function
on $\cone(F)$. These linear functions glue together to
a well--defined concave function $s_{\Pol}\colon V\to \RR$.

The {\it generalized $h$-vector\/} of the polytope $\Pol$ is a
combinatorial invariant defined by  recursion over the faces
of $\Pol$ (see \cite{St2}). In fact, this invariant only depends on the
fan $\Delta_{\Pol}$, and it makes sense to define a generalized
$h$-vector for arbitrary complete fans using the same recursion
formulae. So let us state the definition here in terms of fans.
A {\it fan\/} in a real vector space $V$
is a nonempty set $\Delta$ of strictly convex polyhedral cones  
intersecting pairwise 
in  common faces and such that if a cone belongs to the set $\Delta$
then all its faces also belong to $\Delta$. The fan $\Delta$
is called {\it complete\/} if its support 
$|\Delta|:=\bigcup_{\sigma\in\Delta} \sigma$ equals $V$, and
$\Delta$ is {\it rational\/}
with respect to a lattice $N$ in $\RR$, if all the cones are
generated by vectors in $N$.
For a given
cone $\sigma$, let $V_{\sigma}$ denote the linear span of $\sigma$
in $V$. Let $\Lambda_{\sigma}$ denote the fan  that
we obtain by projecting the boundary of $\sigma$ to  $V_{\sigma}/L$,
where $L$ is a one--dimensional subspace generated by a vector
in the relative interior of $\sigma$.

We introduce two
polynomials, namely $h_{\Delta}$  for each complete fan $\Delta$
and  $g_{\sigma}$ for each strictly
convex polyhedral cone $\sigma$, satisfying
the following recursion:
\begin{enumerate}
\item $g_{0} \equiv 1$
\item $h_{\Delta}(x) =\sum_{\sigma\in\Delta} (x-1)^{\dim \Delta-\dim \sigma} 
g_{\sigma}(x)$ 
\item $g_{\sigma}(x) =\tau_{< [(\dim \sigma)/2]}((1-x) 
h_{\Lambda_{\sigma}}(x)),$
\end{enumerate}
where $\tau_{\leq r}$ denotes the truncation operator $\tau_{\leq r}(\sum_{i=0}^n a_i x^i):=\sum_{i=0}^r a_i x^i$.
The vector formed by the coefficients of the polynomial $h_{\Delta}$  
is called the {\it generalized $h$-vector\/} of the fan $\Delta$.
For a polytope $P$, we set $h_P:=h_{\Delta_P}$. 

For example  $h_n=1$ and $h_{n-1}=|\{\sigma\in\Delta ; \dim\sigma=1\}| -n$,
where $n$ denotes the dimension of the fan.

If $\sigma$ is simplicial cone, then $g_{sigma}=1$. Hence if 
$\Delta=\Delta_P$ is
the fan through the faces of a simplicial polytope $P$, then $h_P=h_{\Delta}=\sum_{\sigma\in\Delta} (x-1)^{\dim \Delta-\dim \sigma}$, which
is nothing but the usual $h$-vector. But for a general polytope, 
its generalized $h$-vector does not only depend on the face numbers but is more
complicated. 

Stanley  showed that the generalized $h$-vector of an $n$-dimensional polytope is palindromic, 
$h_j=h_{n-j}$ for all $j$.
To give an example, the $h$-vector of the $3$-dimensional cross-polytope is 
$(1,3,3,1)$ and
the $h$-vector of its dual, the $3$-dimensional cube is $(1,5,5,1)$.

\section{Combinatorial Intersection Cohomology}

For later use, in this section we very briefly summarize the main results on
the combinatorial intersection cohomology  for fans that are presented
 in \cite{virtual}, and at the same time
we fix the notation. Let $\Delta$ be a 
(not necessarily rational) fan in
a real vectorspace $V=\RR^n$.  

If a subset of cones $\Lambda\subset\Delta$ is again a fan,
then we speak of a {\it subfan\/} and write $\Lambda\prec \Delta$.
In the rational case, where $\Delta$ defines a toric variety,
the subfans of $\Delta$ correspond to the open invariant subsets
of the toric variety, so they define a ``$T$--stable topology''.
That is the motivation for considering the
set of all subfans of an arbitrary fan $\Delta$ together with the empty set
as the open sets of a topology, namely the {\it fan topology\/} on
$\Delta$. A basis of this topology is formed by the {\it affine\/} subfans,
i.e. the subfans that are fans of faces of single cones.
For a cone $\sigma\in\Delta$, we denote the fan of faces of $\sigma$
 by $\langle\sigma\rangle$ and its boundary fan by $\partial\sigma$.

Let $A^{\bullet}:=S^{\bullet}(V^*)$ denote the algebra of 
real--valued polynomial functions on
$V$, together with the grading defined by associating to each linear
function the degree $2$. 
The algebra $A^{\bullet}$ defines a {\it sheaf of graded algebras\/}
$\A$ on $\Delta$ (with the fan topology), where for  $\sigma\in\Delta$
the algebra
$\A(\langle\sigma\rangle)=:A^{\bullet}_{\sigma}$ consists of the
elements of $S^{\bullet}(V_{\sigma}^*)$ viewed as polynomial functions
on $\sigma$. The restriction homomorphisms of $\A$ are given by restriction of
polynomial functions. For $\Lambda\prec\Delta$, 
the sections in $\A(\Lambda)$
correspond to those functions on $\Lambda$ that are conewise
polynomial. Instead of $\A(\Lambda)$ we also write $A^{\bullet}_{\Lambda}$.

Now consider a sheaf $\E$ of graded $\A$--modules on $\Delta$.
To denote the sections $\E(\Lambda)$ of $\E$ on $\Lambda\prec\Delta$ we also
write $E_{\Lambda}$, and we abbreviate $E_{\langle\sigma\rangle}$ to
$E_{\sigma}$. Let $\m$ denote the unique homogeneous maximal
ideal of $\A$. Then forming residue classes modulo $\m$ we
obtain a sheaf of graded real vector spaces $\b{{\cal E}}^{\bullet}$ on $\Delta$,
where $\b{{\cal E}}^{\bullet}(\Lambda):=\b{E}_{\Lambda}$.

The sheaf $\E$ is called a {\it minimal extension sheaf\/} if
 the following properties hold:
\begin{enumerate}
\item $E^{\bullet}_0\simeq \RR^{\bullet}$, where $\RR^{\bullet}$ denotes
$\RR$ viewed as a graded algebra with trivial zero grading.
\item For every $\sigma\in\Delta$, the module $E_{\sigma}^{\bullet}$ is
free over $A_{\sigma}^{\bullet}$.
\item For each cone $\sigma\in\Delta\setminus\{0\}$, the restriction
map $\rho_{\sigma}\colon E_{\sigma}^{\bullet}\to 
E_{\randsigma}^{\bullet}$ induces
an isomorphism $$\b{\rho}_{\sigma}\colon \b{E}_{\sigma}^{\bullet}
\to \b{E}_{\randsigma}^{\bullet}$$
of graded real vector spaces.
\end{enumerate} 

In \cite{virtual}, the authors prove that for any given fan $\Delta$,
 a minimal extension sheaf  exists and is unique
up to isomorphism. If the fan is rational, then we have an associated  toric 
variety $X_{\Delta}$. The equivariant intersection cohomology of
open subsets of $X_{\Delta}$
defines a minimal extension sheaf on $\Delta$ by the assignment
$\Lambda\mapsto \IH_T^*(X_{\Lambda};\RR)$ for $\Lambda\prec\Delta$,
so in particular ${E}^{\bullet}_{\Delta}\simeq \IH^*_T(X_{\Delta};\RR)$.
Moreover,
$\b{E}^{\bullet}_{\Delta}$ is isomorphic to
the intersection cohomology $\IH^*(X_{\Delta};\RR)$ of $X_{\Delta}$ 
(see \cite{Festschrift}).  

If the fan $\Delta$ is not rational, then there is no way of associating
a toric variety to it. But the construction of the sheafs $\E$ and $\b{E}$ still
makes sense, and the authors call it the combinatorial intersection cohomology
of $\Delta$.

Now let us assume that $\Delta=\Delta_{\Pol}$ arises from a polytope {\Pol},
and let $s_{\Pol}$ denote the corresponding strictly concave support
function on $\Delta$. Then the following holds (see \cite{virtual}:

\begin{theorem} ({\rm Barthel, Brasselet, Fieseler, Kaup}) \label{Freiheit} 
 $E^{\bullet}_{\Delta}$ is a free $A^{\bullet}$--module, and
therefore $E_{\Delta}^{\bullet}=
A^{\bullet}\otimes_{\RR} \b{E}_{\Delta}^{\bullet}$.
\end{theorem}

A weak version of the Hard Lefschetz Theorem asserts the
following (see \cite{Karu}):

\begin{theorem} ({\rm Karu}) \label{combHL}
The map $\b{\mu}^{2q}\colon 
\b{E}_{\Delta}^{2q}\to \b{E}_{\Delta}^{2q+2}$
induced by the multiplication with $s_{\Pol}\in {\cal A}^2(\Delta)$
is injective for $2q\leq n-1$ and surjective for $2q\geq n-1$.
\end{theorem}

In \cite{virtual} the authors showed that if \ref{combHL} is true then
the Betti-numbers of the combinatorial intersection cohomology are
precisely given by the coefficients of the $h$-polynomial of the fan:

\begin{theorem}({\rm Barthel, Brasselet, Fieseler, Kaup}) \label{Bettizahlen}
$$h_{\Delta}(t^2)=\sum_{q=0}^{2n} (\dim \b{E}_{\Delta}^{q})\, t^q\,.$$
\end{theorem}

The fact that the $h$-polynomial is palindromic is reflected by a combinatorial
version of Poincar\'e-duality, also proved by the four authors.

For $\E$ and $\b{{\cal E}}^{\bullet}$ 
we can define a Poincar\'e--series $\IQ_{\Delta}$
and a Poincar\'e--polynomial $\IP_{\Delta}$ respectively as follows:
$$\IQ_{\Delta}(t):=\sum_{q\geq 0} (\dim E_{\Delta}^{q}) \,t^{q}
\quad\text{and}\quad
\IP_{\Delta}(t):=\sum_{q\geq 0} (\dim \b{E}_{\Delta}^{q})\, t^{q}\,.$$
Note that here we do not follow the convention used in \cite{virtual}
in order to be consistent with the notation used in \cite{hVektoren}.
One obtains the Poincar\'e--series used in \cite{virtual} from
ours by viewing it as a function in $t^2$.   

As a consequence
of the first part of the above theorem we obtain:
$$ \IQ_{\Delta}(t)={1\over (1-t^2)^n}\mal \IP_{\Delta}(t)\,.\leqno{(1)}$$

\section{Refined Poincar\'e--Series}

From now on let $\Delta$ denote a  complete fan, and
assume that for every $\sigma\in\Delta$ also $-\sigma\in\Delta$,
in other words assume that $\Delta$ is {\it centrally symmetric\/}.
Let $\varphi=-\id_V$ denote the central symmetry.
Being an invertible linear transformation,  $\varphi$ 
induces an $\RR$--linear   automorphism of the graded algebra
$A^{\bullet}=S^{\bullet}(V^*)$. Note that for every cone $\sigma\in\Delta$,
we have $V_{\sigma}=V_{-\sigma}$. Since $A^{\bullet}_{\sigma}$
is the algebra of polynomial functions on $V_{\sigma}$ restricted to $\sigma$,
the algebras $A^{\bullet}_{\sigma}$ and $A^{\bullet}_{-\sigma}$ are
not identical, but canonically isomorphic.
The action of $\varphi$ on $V_{\sigma}$ induces
an isomorphism of graded algebras from $A^{\bullet}_{\sigma}$ to
$A^{\bullet}_{-\sigma}$ that is compatible with this canonical
isomorphism. Moreover, the induced isomorphisms are 
compatible with the restriction homomorphisms
$\rho^{\sigma}_{\tau}\colon A^{\bullet}_{\sigma} \to A^{\bullet}_{\tau}$
for every $\tau\prec\sigma$.
That means that in fact $\varphi$ defines a natural automorphism
of $\Delta$ as a ringed space equipped with the
sheaf of graded algebras $\A$. We can also define an action
of $\varphi$ on the minimal extension sheaf $\E$ on $\Delta$.

\begin{lemma}
There are isomorphisms of graded vector spaces
$$\varphi\colon E^{\bullet}_{\sigma}\to E^{\bullet}_{-\sigma}$$
for every $\sigma\in\Delta$ that are equivariant with respect
to the module structure
over  $A^{\bullet}_{\sigma}$ and 
$A^{\bullet}_{-\sigma}$ respectively and
compatible with the restriction homomorphisms of the sheaf $\E$.
\end{lemma}

\proof To define the required isomorphisms, we proceed by recursion over
the $k$--skeleton $\Delta^{\leq k}$ of $\Delta$
following the recursive construction of $\E$ as
in Section~1 of \cite{virtual}. On 
$E^{\bullet}_0=\RR^{\bullet}$,
where $0$ denotes the zero cone, $\varphi$ acts as the identity.
Now assume that the isomorphisms have been defined for $\Delta^{<k}$,
and let $\sigma\in\Delta^k$. We can assume that 
$E_{\sigma}^{\bullet}=A^{\bullet}_{\sigma}
\otimes_{\RR} \b{E}^{\bullet}_{\randsigma}$.
By induction, we already have an isomorphism
$\varphi\colon E^{\bullet}_{\randsigma}\to E^{\bullet}_{-\randsigma}$,
and since the maximal ideal $\m$ of $\A$ is $\varphi$--stable,
$\varphi$ induces an isomorphism
$\b{\varphi}\colon \b{E}^{\bullet}_{\randsigma}
\to \b{E}^{\bullet}_{-\randsigma}$. Together with
the map from $A^{\bullet}_{\sigma}$ to
$A^{\bullet}_{-\sigma}$ determined by $\varphi$,
that provides us with an isomorphism of graded vector spaces
$\varphi\colon E^{\bullet}_{\sigma}\to E^{\bullet}_{-\sigma}$.
By construction, this map is equivariant as a map
from an $A^{\bullet}_{\sigma}$--module to an 
$A^{\bullet}_{-\sigma}$--module. 

In the construction of $\E$,
the restriction homomorphism from $E^{\bullet}_{\sigma}=
A^{\bullet}_{\sigma}
\otimes_{\RR} \b{E}^{\bullet}_{\randsigma}$
to $E^{\bullet}_{\randsigma}$ is defined using the restriction homomorphism
$\rho^{\sigma}_{\randsigma}\colon A^{\bullet}_{\sigma} \to
A^{\bullet}_{\randsigma}$ and
an $\RR$--linear section $s_{\sigma}\colon \b{E}^{\bullet}_{\randsigma} \to
E^{\bullet}_{\randsigma}$ of the residue class map
$E^{\bullet}_{\randsigma} \to \b{E}^{\bullet}_{\randsigma}$.
The section $s_{\sigma}$ can be chosen freely. So we can assume
that for any pair $\sigma,-\sigma$ of antipodal cones in $\Delta$
the corresponding sections have been chosen such that the
following diagram is commutative:
$$
\matrix{%
\b{E}^{\bullet}_{\randsigma} & \bigtopmapright{s_{\sigma}} &
E^{\bullet}_{\randsigma} &\cr
\lmapdown{\varphi} & & \rmapdown{\varphi} &\cr
\b{E}^{\bullet}_{-\randsigma}& \bigtopmapright{s_{-\sigma}} & 
E^{\bullet}_{-\randsigma} &. \cr}
$$
That implies 
compatibility of $\varphi$ with the restriction homomorphisms of $\E$.
\endproof

In particular, we obtain an induced automorphism $\varphi$ on the
module $E_{\Delta}^{\bullet}$ of global sections of $\E$,
and though the automorphism is not canonical, 
 the dimensions of the eigenspaces in each graded piece are
uniquely determined and therefore the so--called
{\it refined Poincar\'e--series\/} for the action of $\varphi$
on $\E$ depends only on $\Delta$.
This series is defined as a polynomial over
the group ring $\ZZ[G]$ of the character group $G:=\{\pm 1\}$ 
 of the
group generated by $\varphi$ in $\GL(V)$, namely:
$$\IQ_{\Delta}^{\varphi}(t):=
\sum_{q\geq 0} (\dim (E_{\Delta}^{q})^+  + 
\dim (E_{\Delta}^{q})^- \chi) \,\, t^{q}\,,$$
where $\chi$ denotes the element corresponding to $-1$ in
$\ZZ[G]$, and
the superscripts $+$, $-$ indicate the eigenspaces for the eigenvalues
$+1$ and $-1$ respectively. The refined Poincar\'e--polynomial
$\IP_{\Delta}^{\varphi}$ for the action of $\varphi$ on 
$\b{{\cal E}}^{\bullet}$
is defined analogously.

By Theorem \ref{Freiheit}, $E_{\Delta}^{\bullet}$
is a free $A^{\bullet}$--module. 
The action of $\varphi$ on 
$\b{E}_{\Delta}^{\bullet}$ is induced by taking residue classes.
So if we choose a homogeneous basis for $\b{E}_{\Delta}^{\bullet}$
and preimages under the residue class map in $E_{\Delta}^{\bullet}$
to define an isomorphism
$$E_{\Delta}^{\bullet}\to \b{E}_{\Delta}^{\bullet}\otimes_{\RR} A^{\bullet}\,,
$$
then this isomorphism is automatically compatible with the
action of $\varphi$.
That implies
$$ \IQ_{\Delta}^{\varphi}(t)=
{1\over (1-\chi t^2)^n}\mal \IP_{\Delta}^{\varphi}(t)\,.
\leqno{(2)}$$

To obtain a relation between the  Poincar\'e--series $\IQ_{\Delta}$
and its refined version $\IQ_{\Delta}^{\varphi}$, we can use the fact that 
the minimal extension sheaf $\E$
as a sheaf of real vector spaces can be written
 as a direct sum of simpler subsheafs.
Note that $\E$ is a {\it flabby\/} sheaf on $\Delta$. Here that means
that the restriction homomorphism $\rho^{\sigma}_{\randsigma}\colon
E^{\bullet}_{\sigma} \to E^{\bullet}_{\randsigma}$ is surjective
for all $\sigma\in\Delta$.

For $\sigma\in\Delta$, let $\I_{\sigma}$ denote the
{\it characteristic sheaf\/} of $\sigma$ defined on
 $\Lambda\prec\Delta$ by
$$\I_{\sigma}(\Lambda):=\cases{\RR & if $\sigma\in\Lambda$\cr
                                      0 & otherwise} \,.$$
Then there is an isomorphism of sheafs of graded real vector spaces
$$\E\simeq \bigoplus_{\sigma\in\Delta} \I_{\sigma}\otimes_{\RR}
K_{\sigma} \,,
\leqno{(3)}$$
where $K_{\sigma}$ denotes the kernel of the restriction homomorphism
$\rho^{\sigma}_{\randsigma}\colon E^{\bullet}_{\sigma}\to
E^{\bullet}_{\randsigma}$ (see Section~3, \cite{virtual}).
 
For every $\sigma\in\Delta$, $\varphi$ induces a map from
$\I_{\sigma}$ to $\I_{-\sigma}$ that is compatible with the
action of $\varphi$ on $\Delta$. And using these maps together
with the  action of $\varphi$ on $\E$, we obtain an induced
$\varphi$--action on the direct sum on the
righthandside of (3),
where $\varphi$ maps $\I_{\sigma}$ to $\I_{-\sigma}$ and
$K_{\sigma}$ to $K_{-\sigma}$. Modifying the proof of the decomposition
theorem from
\cite{virtual}, we can show the following:

\begin{lemma} The isomorphism
of sheafs of graded real vector spaces on $\Delta$ in (3)
can be chosen as $\varphi$--equivariant.
\end{lemma}

\proof
We prove our claim by induction on the number of antipodal pairs of
cones in $\Delta$. Suppose that there is a $\varphi$--equivariant
decomposition of $\E$ into $\varphi$--stable flabby sheafs
$$\E\simeq\F\oplus \left(\bigoplus_{\sigma\in \Lambda} \I_{\sigma}\otimes_{\RR}
K_{\sigma}\right)\,,$$
where the sum is taken over a $\varphi$--stable subset $\Lambda$ of $\Delta$
(that is not necessarily a subfan),
such that $\F(\sigma)=0$ for all $\sigma\in\Lambda$.

Then choose a pair of antipodal cones $\sigma,-\sigma
 \in \Delta\setminus\Lambda$
of minimal dimension $k$ with $\F(\sigma)\ne 0\ne \F(-\sigma)$. 
We  have to show that we can write  $\F$ 
 as a direct sum of $\varphi$--stable
flabby subsheafs
$\F=\G\oplus \H$ such that $\H(\sigma)=0$ and if $\sigma=0$ then
 $\G\simeq\I_0 \otimes_{\RR} K_0$ and if $\sigma\ne 0$
then $\G\simeq (\I_{\sigma}\otimes_{\RR}
K_{\sigma}) \oplus (\I_{-\sigma}\otimes_{\RR}
K_{-\sigma})$. We define $\G$ and $\H$
on the $k$--skeleton $\Delta^{\leq k}$ by
$${\G}(\tau):=\cases{K_{\tau} & if $\tau=\pm\sigma$\cr
                      0          & otherwise}
\quad\text{and}\quad
{\H}(\tau):=\cases{0 & if $\tau=\pm\sigma$\cr
                                      \F(\tau) & otherwise}\,. $$
Now suppose, that $\F=\G\oplus\H$ is already defined on $\Delta^{\leq m}$,
and consider a pair of antipodal cones $\pm \tau$ of dimension
$m+1$. If $\sigma$ is neither a face of $\tau$ nor of $-\tau$,
then set $\G(\pm\tau)=0$ and $\H(\pm\tau):=\F(\pm\tau)$.
Otherwise say $\sigma\prec\tau$ and $-\sigma\prec-\tau$.
Note that $\tau$  cannot contain both $\sigma$ and $-\sigma$
as a face. 

By assumption, we have a decomposition 
$\F(\randtau)=\G(\randtau)\oplus \H(\randtau)$, and
$\G(-\randtau)=\varphi(\G(\randtau))$ and
$\H(-\randtau)=\varphi(\H(\randtau))$.
Since $\F$ is flabby, we can choose a decomposition
$\F(\tau)=U\oplus W$ such that the restriction homomorphism
$\rho^{\tau}_{\randtau}\colon \F(\tau)\to \F(\randtau)$ induces
an isomorphism $U\to \G(\randtau)$
and a surjective homomorphism from $W$ to $\H(\randtau)$.
Since the action of $\varphi$ is compatible with the restriction
homorphisms,
for the decomposition $\F(-\tau)=\varphi(U)\oplus\varphi(W)$
the following holds: The restriction homomorphism
$\rho^{-\tau}_{-\randtau}\colon \F(-\tau)\to \F(-\randtau)$ induces
an isomorphism $\varphi(U)\to \G(-\randtau)$
and a surjective homomorphism from $\varphi(W)$ to $\H(-\randtau)$.
Now set $\G(\tau)=U$ and $\G(-\tau)=\varphi(U)$,
$\H(\tau)=W$ and $\H(-\tau)=\varphi(W)$. Then $\G$ and $\H$ have
the required properties.
\endproof

Now consider the action of $\varphi$ on the direct sum on the
righthandside of (3).
Apparently, for every cone $\sigma\ne 0$, the action of $\varphi$
interchanges the summands $\I_{\sigma}\otimes_{\RR}
K_{\sigma}$ and $\I_{-\sigma}\otimes_{\RR}
K_{-\sigma}$. So for every $q$, in 
$\bigoplus_{\sigma\ne 0} \I_{\sigma}\otimes_{\RR}
K_{\sigma}$
the eigenvalue $+1$ and $-1$ occur with the same multiplicity.
We obtain the relation
$$\IQ_{\Delta}^{\varphi}(t)-1=
{1\over 2}((\IQ_{\Delta}(t)-1)+(\IQ_{\Delta}(t)-1)\mal \chi) =
{1\over 2} (1+\chi) (\IQ_{\Delta}(t)-1) \,.
\leqno{(4)}$$

\section{Lower Bounds for the Generalized $h$--Vector}

Summarizing the considerations in the previous section, we obtain
the following description of the refined Poincar\'e--polynomial:

\begin{proposition}
Let $\Delta$ be a centrally symmetric complete fan of dimension $n$. Then
$$\IP_{\Delta}^{\varphi}(t)= {1\over 2}(\IP_{\Delta}(t)+(1+t^2)^n)
+ {1\over 2} \chi (\IP_{\Delta}(t)-(1+t^2)^n)\,.$$
\end{proposition}

\proof
Inserting (1) in (4), we obtain
$$\IQ_{\Delta}^{\varphi}(t)={1\over 2}(1+\chi)\IQ_{\Delta}(t) + 
{1-\chi \over 2}\,.$$
Using (2), that implies
$$\IP_{\Delta}^{\varphi}(t)= 
{1\over 2}{(1+\chi)(1-\chi t^2)^n \over (1-t^2)^n} \IP_{\Delta}(t) 
+  {1-\chi \over 2}(1-\chi t^2)^n\,.
$$
Note that since $\chi^2=1$, we have $(1-\chi t^2)(1+\chi)=(1-t^2)(1+\chi)$
and $(1-\chi t^2)(1-\chi)=(1+t^2)(1-\chi)$. This implies
$$\IP_{\Delta}^{\varphi}(t)={1-\chi \over 2}(1+t^2)^n +
{1+\chi \over 2} \IP_{\Delta}(t)\,. \quad\kasten$$

\bigskip

We now  apply this proposition 
to polytopal centrally symmetric fans.


\begin{theorem}\label{Hauptsatz}
Let ${\Pol}$ be a centrally--symmetric polytope of dimension $n$. 
Then the polynomial
$$h_{\Pol}(x)-(1+x)^n\,,$$
has nonnegative, even coefficients, it is palindromic and unimodal.
In particular, we have
the following bounds for the coefficients
$h_j$ of $h_{\Pol}$:
$$h_j-h_{j-1}\geq {n\choose j}- {n\choose j-1}\quad\text{for $j=1,\dots,n/2$.}
$$
\end{theorem}

\proof
As before let $\Delta:=\Delta_{\Pol}$ denote the fan through the 
faces of ${\Pol}$ and
let $s_{\Pol}$ denote the $\Delta$--strictly convex support function defined
by ${\Pol}$. 
By Theorem \ref{Bettizahlen}, $\IP_{\Delta}(t)=h_{\Pol}(t^2)$.
The symmetry follows immediately from the combinatorial Poincar\'e--duality.
Moreover, by the above proposition, we have
$${1\over 2}(\IP_{\Delta}(t)-(1+t^2)^n)=
\sum_{q\geq 0} (\dim (\b{E}_{\Delta}^{q})^-)\, t^{q}\,.$$
This implies in particular, that all the coefficients of
$p(t):=\IP_{\Delta}(t)-(1+t^2)^n$ are nonnegative and even. 

Since the support
function $s_{\Pol}\in {\cal A}^2(\Delta)$ is invariant under $\varphi$,
we have $$s_{\Pol}\cdot (\b{E}_{\Delta}^{q})^- \subset \b{E}_{\Delta}^{q+2})^-\,.$$
Now it follows from the combinatorial Hard Lefschetz theorem 
(see Theorem \ref{combHL}) that $\dim ((\b{E}_{\Delta}^{2q})^-)\leq
\dim ((\b{E}_{\Delta}^{2q+2})^-)$ for $2q\leq n-1$, and that means
that the polynomial $p$ is unimodal.
\endproof

Note that $(1+x)^n$ occurs as the $h$-polynomial of the $n$-dimensional cross-polytope. We can reformulate the lower bounds given by the theorem in terms of the partial ordering on real polynomials of degree $n$ defined by  coefficientwise comparison, i.e. $a=\sum_{j=0}^n a_j x^j \leq b=\sum_{j=0}^n b_j x^j$ if and only if $a_j\leq b_j$ for all $j$.  
The $h$-polynomial of the $n$-dimensional cross-polytope is minimal in this sense and in fact this is the only polytope realizing the minimum.

\begin{corollary}\label{Minimum}
Let $\Pol$ be an $n$-dimensional centrally symmetric polytope. Then
$$h_{\Pol}\geq (1+x)^n\,.$$ 
Moreover, equality holds if and only if the polytope $\Pol$ is affinely equivalent to the $n$-dimensional cross-polytope.
\end{corollary}

\proof Suppose that $h_{\Pol}=(1+x)^n$ for a centrally symmetric polytope $\Pol$. Then in particular, $h_{n-1}=n$.
On the other hand, as mentioned in Section~1, we always
have, $h_{n-1}=
|\{\hbox{\rm vertices of\ } \Pol\}|-n$. So $\Pol$ has $2n$ vertices. Let us
choose a facet $F$ of $\Pol$. Since $F$ contains at least $n$ vertices,
and it is disjoint from its opposite facet $-F$, we obtain that $\Pol$ is
the convex hull of $F$ and $-F$. That means that $\Pol$ is affinely equivalent
to the $n$-dimensional cross-polytope.  
\endproof

\bibliography{}

\begin{thebibliography}{KaStZe}%
%
\bibitem[AC]{hVektoren} A.~A'Campo--Neuen,
{\it On Generalized $h$--Vectors of Rational Polytopes with
a Symmetry of Prime Order\/}, Discrete Comput. Geom. {\bf 22} (1999), 259--268.
%
%
\bibitem[Ad]{Ad} R.~M.~Adin, {\it On face numbers of rational simplicial polytopes with symmetry\/}. Adv. Math.~{\bf 115} (1995), 269--285.

\bibitem[BBFK1]{Festschrift} G.~Barthel, J.-P.~Brasselet, K.-H.~Fieseler, L.~Kaup,
{\it Equivariant intersection cohomology of toric varieties\/}. In: Algebraic Geometry:
Hirzebruch 70. Contemp. Math. AMS {\bf 241} (1999), 45--68.

\bibitem[BBFK2]{virtual} G.~Barthel, J.-P.~Brasselet, K.-H.~Fieseler, L.~Kaup,
{\it Combinatorial intersection cohomology for fans\/}. T\^ohoku Math. J. 54 (2002), 1--41.

\bibitem[BL]{Bressler} P.~Bressler, V. A. Lunts, {\it Intersection Cohomology on Nonrational Polytopes\/}, preprint math.AG/0002006.

\bibitem[Fi1]{Fi} K.-H.~Fieseler, {\it Rational intersection cohomology of 
projective toric varieties\/}. J. reine angew.~Math.~{\bf 413} (1991), 88--98.

\bibitem[Ka]{Karu} K. Karu, {\it Hard Lefschetz Theorem for nonrational polytopes\/}.
preprint math.AG/0112087 v3.

\bibitem[St1]{St1} R.~Stanley, {\it On the number of faces of centrally--symmetric simplicial polytopes}, Graphs Combin.~{\bf 3} (1987), 55--66.

\bibitem[St2]{St2}  R.~Stanley, {\it Generalized $h$-vectors, intersection cohomology of toric varieties, and related results}, in: Commutative Algebra and Combinatorics (M.~Nagata and H. Matsumura, eds.), Advanced Studies in Pure Math..~{\bf 11}, Kinokuniya, Tokyo, and North--Holland, Amsterdam/New York, 1987, pp.~187--213. 
%
\end{thebibliography}

\end{document}